\documentclass[nofootinbib,
10pt,
 amsmath,amssymb,
 aps, physrev,
]{revtex4-2}

\usepackage{graphicx}
\usepackage{dcolumn}
\usepackage{bm}
\usepackage{hyperref}

\usepackage{braket} 
\usepackage{pythonhighlight} 
\usepackage{algorithm2e}

\usepackage[textsize=footnotesize,backgroundcolor=red!25]{todonotes}

\newcommand{\mono}[1]{\texttt{#1}}
\renewcommand{\vec}[1]{\boldsymbol{#1}}
\newcommand{\fcaption}[1]{\caption{\small #1}}

\begin{document}

\title{\textbf{Numerical methods for the simulation of quantum walks and quantum annealing} 
}%

\author{Asa Hopkins}
\email{Contact author: asa.hopkins@strath.ac.uk}
\author{Viv Kendon}%
\affiliation{%
Department of Physics and SUPA, University of Strathclyde, Glasgow G4 0NG, United Kingdom
}%

\date{\today}%

\begin{abstract}

It is known that Chebyshev based polynomial approximation gives a near-optimal rate of convergence for calculating a function of a Hermitian matrix.  However, previous work has not discussed the option of true minimax approximation, nor the specifics of writing a high performance implementation with a rigorous analysis of errors. This work provides such an analysis and an open-source implementation of three approximation methods in C++.

\end{abstract}

\maketitle

\section{Introduction and Background}

For general matrices, a matrix exponential can be calculated by careful evaluation of the Taylor series, as done in \cite{Higham}.  For common special cases, more efficient methods can significantly reduce the computational costs.   A motivating example will be taken from quantum mechanics. The evolution of a quantum system is governed by the time-dependent Shr\"odinger equation, 
\begin{align}\label{eq:schroed}
    i \hbar \frac{d}{dt}\ket{\psi}= \hat{H}(t) \ket{\psi},
\end{align}
where $i = \sqrt{-1}$, $t$ is time, $\ket{\psi}$ is a vector representing the quantum state, and $\hat{H}$ is a matrix representing the system Hamiltonian. For an $n$ qubit system, $\ket{\psi}$ is of length $N = 2^n$ (qubits are two-state quantum systems).
For a time-independent Hamiltonian, this can be solved exactly to give
\begin{align}\label{eq:SchroedSol}
    \ket{\psi(t)}= \exp(-it \hat{H}) \ket{\psi}
\end{align}
in units where $\hbar = 1$.
In this work, we consider Hamiltonians used for quantum optimization: $\hat{H}$ will be chosen from Hamiltonians of the form $\hat{H} = \hat{H}_P + \gamma \hat{H}_G$, where $\hat{H}_G$ is defined using Pauli-X gates $\hat{X}_j$ on the $j$th qubit as 
$$\hat{H}_G = \sum_{j=0}^{n-1} \hat{X}_j,$$ 
and $\hat{H}_P$ is a diagonal matrix. This Hamiltonian is sparse, with only $n+1$ non-zero elements per row.  Hamiltonians with this level of sparsity are common in quantum computing, due to the physical restrictions on interactions between qubits in quantum hardware.
The goal of quantum optimization is to perform a quantum evolution which ends with $\ket{\psi}$ in a low energy state of $\hat{H}_P$ -- for more details, see the companion paper \cite{hopkins}. 
For this work, the focus will be on how to best evaluate $\exp(-it\hat{H})$ classically for some sparse Hermitian $\hat{H}$.

\begin{figure}
\includegraphics[width=0.60\linewidth]{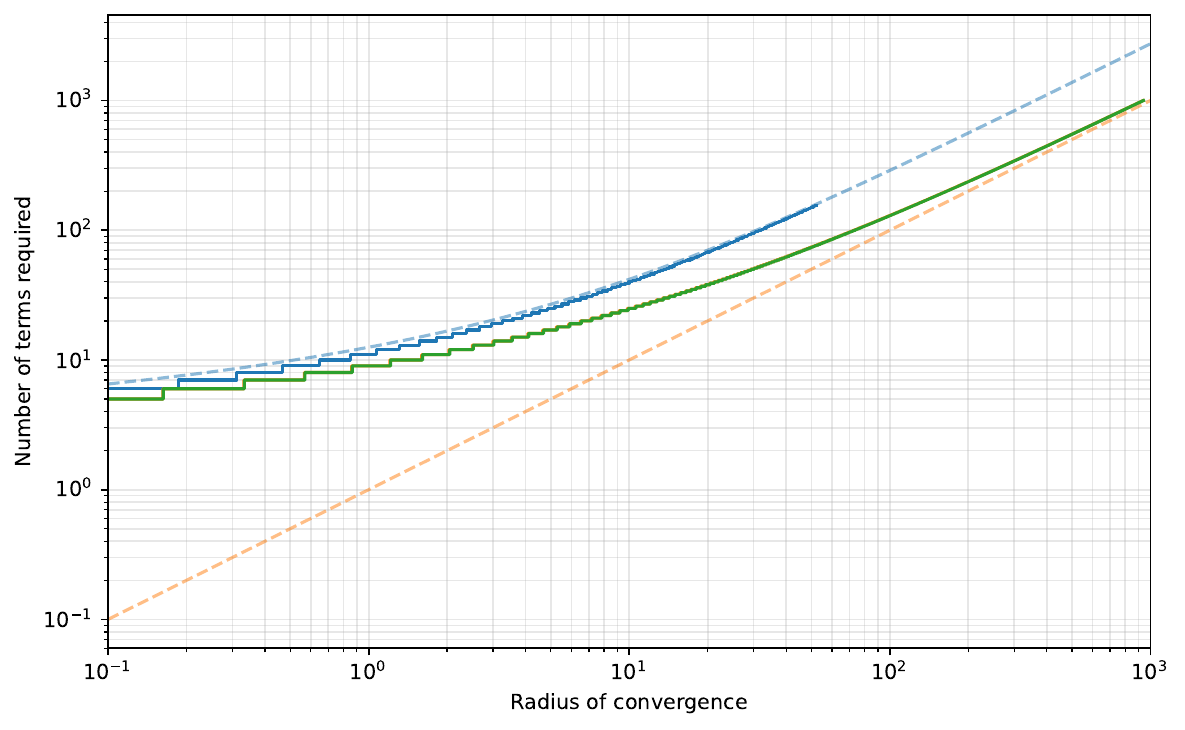}
\fcaption{The number of terms required for different polynomial expansions of the function $\sin(x) + \cos(x)$ to reach a truncation error of $2^{-23}$ corresponding to 32-bit floating point precision, as a function of the required radius of convergence. The solid blue line is the Taylor series, the solid green line is the minimax polynomial calculated using the CF method. A solid orange line is drawn for the Chebyshev series but the green line overlaps it closely. The dashed blue line is the asymptote for the number of terms for the Taylor series, and the dashed orange line is the asymptote for the number of terms for the Chebyshev series.}
\label{fig:convergence}
\end{figure}

Numerically solving the ODE in Eq.~\eqref{eq:schroed} for time independent Hamiltonians, using Eq.~\eqref{eq:SchroedSol}, reduces to evaluating the action of a matrix exponential on a vector, which can be done more efficiently than evaluating $\exp(-it\hat{H})$ directly or diagonalizing the matrix. 
In the quantum case, the matrices are guaranteed to be Hermitian, and therefore the eigenvalues are guaranteed to be real. If a bound can be put on the eigenvalues, then a Chebyshev series can be applied which converges in much fewer terms than a Taylor series, as shown in Fig.~\ref{fig:convergence}. 
When the matrix exponential itself is of interest, \cite{Higham2} discusses other techniques that can converge in fewer terms, such as scaling and squaring, or Pad\'e approximations, but these are only competitive when the matrix exponential itself is needed, rather than just its action on a vector.

While use of Chebyshev polynomials in this context has been known and discussed since the mid-1980s \cite{propagate,spectral,1989}, efficient implementations are not yet widely available in open-source software packages.
The contributions of our work are: rigorous analytic error bounds on Chebyshev approximation for matrices via the Clenshaw \cite{Clenshaw} algorithm, true minimax approximation found with the Carath\'eodory-Fej\'er algorithm, and a concrete implementation explaining the optimisations techniques found. Our implementation is specialised for quantum walk evolution, currently supporting Ising models on hypercube graphs. The Clenshaw algorithm optimisation and the matrix-free hypercube optimisation are both of wider interest and applicability.
The code used in this paper can be found at \cite{AsaCode}.

The remainder of the paper is organized as follows.  In section \ref{ssec:prior} we discuss prior work on this and related algorithms.  In section \ref{sec:methods} we provide details of our numerical and analytical methods, and the resulting improvements in runtimes and convergence are illustrated numerically.  In section \ref{sec:diss} we summarise and discuss directions for further work and implementation.

\subsection{Prior work}\label{ssec:prior}

\citeauthor{propagate} \cite{propagate} lay the groundwork for the Chebyshev method. 
\citeauthor{spectral} \cite{spectral} discusses the need for Clenshaw-based evaluation in the Chebyshev basis, to avoid catastrophic cancellation and discusses techniques for memory reduction.
However, they only measure errors empirically and do not provide a modern forward or backward error analysis.  They also only cover the case of scalar evaluations with the Clenshaw algorithm: adapting to matrix-vector multiplications is not trivial since there are optimizations that are only applicable in the matrix-vector case. 
The Chebyshev interpolant has been proven to be within a factor $O(\ln k)$ of the error of the true minimax polynomial  \cite{lebesgue}, where $k$ is the degree of the polynomial,
and \cite{spectral} states that calculation of the true minimax polynomial is too expensive to be worthwhile. 
This is true in the case of solving scalar PDEs, but for other cases, such as the quantum mechanical example discussed here, the effort to find the true minimax polynomial can be small compared to even a single matrix-vector multiplication.

\citeauthor{Laguerre} \cite{Laguerre} focus on using Laguerre and Hermite polynomials to calculate the action as they are valid on unbounded regions and therefore avoid needing a bounded region containing the eigenvalues. This gives a Clenshaw-like algorithm and a \emph{filtered conjugate residual} (FCR) algorithm for evaluating matrix functions in general orthogonal polynomial bases. They conclude that the FCR algorithm has fewer floating point operations.  However, it is memory accesses that are important, and their implementation of the Clenshaw algorithm is not optimal for memory accesses. Even so, it may be that an optimised FCR method would be competitive to the Clenshaw algorithm.

\citeauthor{kernel} \cite{kernel} also discuss Chebyshev methods in the context of approximating kernel functions, and discuss techniques for improving on the Chebyshev series when the domain contains a discontinuity. The quantum mechanical example we use is a smooth function which doesn't require such techniques.
\citeauthor{KITE} \cite{KITE} use Chebyshev based methods for density of states calculations and Green's functions of a Hamiltonian, rather than for evolution of a state. Their KITE package fixes some of the optimisation limitations of \cite{propagate,spectral}, with a big focus on memory locality and reducing cache misses. However, they do not provide strict error bounds or explicitly discuss the necessary modifications to the Clenshaw algorithm.

\citeauthor{HighamBook2} \cite{HighamBook2} also mentions that for Hermitian matrices, Chebyshev polynomials are best for large sparse matrices, referring to \citeauthor{1989} \cite{1989}, which provides an error analysis of the method. This analysis is largely based on the error analysis of the scalar function approximation and doesn't take floating point error into account. 
Although Chebyshev-based evaluation has long been known in relevant literature \cite{matrices,cheb,cheb2}, it seems to have been overlooked in practice, as it doesn't currently feature in relevant open-source software packages. Software such as QuTiP \cite{QuTiP} uses numerical integration techniques, which are less stable and can fail to converge if given too large a time step, whereas software such as QuEST \cite{quest} uses Trotter-Suzuki \cite{TrotterSuzuki} decomposition instead as it is restricted to unitary methods. Trotter-Suzuki could be used here, but the conclusion in \cite{Trotter} is that Chebyshev expansion is preferred where it is available.

\section{Methods and results} \label{sec:methods}

\subsection{Calculating the action of matrix functions}
To evaluate Eq.~\eqref{eq:SchroedSol}, 
where the vectors and matrices are of dimension $N=2^n$, for small $n$ it is viable to diagonalise $\hat{H}$.  Then one can transform $\ket{\psi}$ into the basis where $\hat{H}$ is diagonal, apply the exponential, and transform back. This requires $O(N^3)$ operations since $\hat{H}$ is Hermitian. This can be done using the \mono{NumPy} \cite{NumPy} function \mono{numpy.linalg.eigh}. In our benchmarks, sparse methods from SciPy \cite{SciPy} become faster beyond around $n = 12$.
The library \mono{CuPy} \cite{CuPy} offers GPU-accelerated versions of many \mono{NumPy} and \mono{SciPy} functions, including a \mono{cupy.linalg.eigh} function. However, in our benchmarks the overhead of GPU offloading means this is not faster than other methods at any value of $n$ tested. GPU implementation can be faster, but care must be taken to mitigate the penalty of moving data between system RAM and video RAM (VRAM). The mitigations that have been tested here are to batch simulations together, and to generate the $\hat{H}_P$ matrices directly in VRAM where possible. This allows making use of the faster RAM available in GPUs, which can speed up computations by a factor of 20 even while being more energy efficient.

To calculate $\exp(-i\hat{H}t) \ket{\psi}$ without explicitly calculating $\exp(-i\hat{H}t)$ -- which would require dense matrix multiplication even if $\hat{H}$ is sparse -- the typical technique is to use Horner's algorithm to evaluate enough terms of the Maclaurin series of $\exp(x)$ to reach the desired accuracy. Let $a$ be an array of length $m$ such that $p(x) = \sum_{j=0}^{k} a_j x^j$ is the desired polynomial, then $p(H) \ket{\psi}$ can be evaluated using Algorithm \ref{alg:horner}.

\begin{algorithm}
\begin{minipage}{0.95\linewidth}
\begin{python} 
def horner(p, H, psi):
    """Evaluate p(H) @ psi using Horner's method.
    Parameters:
    ----------
    p : list[float] of coefficients, where p[m] is the coefficient of x^m
    psi : NDArray[np.floating] vector to which p(H) is applied.
    H : Any linear operator implementing `@`
    (matrix-vector product) and compatible with psi.
    """
    output = p[-1]*psi.copy()
    for coeff in p[-2::-1]:
        output = H @ output + coeff*psi
    return output
\end{python}
\caption{Horner's method for evaluating \mono{p(H) @ psi}, in Python.}
\label{alg:horner}
\end{minipage}
\end{algorithm}

An issue arises since the Maclaurin series of $\exp(x)$ is not numerically stable for large $x$.  To solve this, an integer $s$ is chosen so that the calculation of $\exp(\hat{H}/s) \ket{\psi}$ is stable, and then repeated $s$ times to give the desired calculation. This is the approach taken in \cite{Higham}. It is known that the best polynomial approximation on the unit disc is the Maclaurin polynomial, so no better choice is possible in the most general case \cite{taylor}.
For our quantum mechanical example, however, this is sub-optimal as $\hat{H}$ is Hermitian and has entirely real eigenvalues. This means the polynomial $p$ can be restricted to only approximate $\exp(ix) = \cos(x) + i\sin(x)$ for $x \in \mathbb{R}$, for which better options exist.

\begin{figure}
\includegraphics[width=0.90\linewidth]{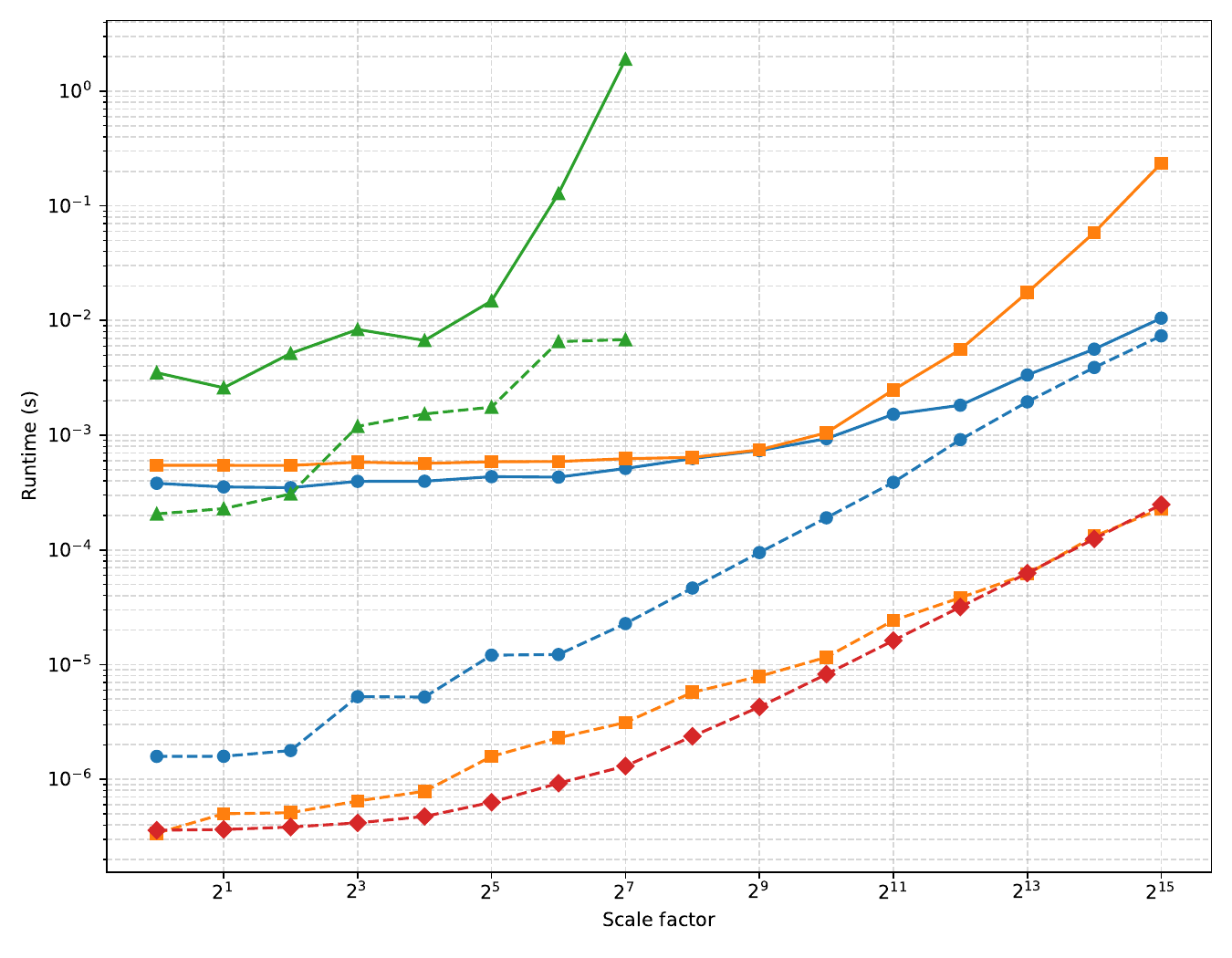}
\fcaption{The runtime needed to create a polynomial $p(x)$ which approximates $\sin(x) + \cos(x)$ on the interval [-scale,scale] to an error of $2^{-23}$ using a variety of methods. Solid lines show Chebfun algorithms and dashed lined show our C++ implementations. Green triangles show the Remez algorithm, Blue circles show interpolation at the Chebyshev nodes, orange squares show the CF method and red diamonds show Miller's algorithm. It can be seen that Miller's algorithm combined with our CF implementation give vastly improved runtime over the original Chebfun methods, making them viable for practical use.}
\label{fig:runtime2}
\end{figure}

The book \textit{Approximation Theory and Approximation Practice} \cite{ATAP} provides a comprehensive survey of approximation techniques, from which a few results will be used. The solution to the problem of finding the best polynomial approximant to a function on a fixed real interval is called the minimax problem, which can be solved using the Remez exchange algorithm. Trefethen argues \cite[Chapter 16]{ATAP} that using the Remez algorithm is too computationally expensive to be worthwhile, and favours Chebyshev interpolation which is a factor of $O(\ln k)$ 
from minimax, where $k$ is the degree of the polynomial, and can be calculated easily with a fast Fourier Transform (FFT) \cite{Fourier}. The transform is equivalent to the Discrete Cosine Transform, which has a more efficient implementation given by \citeauthor{DCT} \cite{DCT}. There is an alternative method for calculating minimax polynomials, called the Carath\'eodory-Fej\'er (CF) method, which is described in \cite{CF} and \cite[Chapter 20]{ATAP}. In simple terms, this method takes a Chebyshev expansion of length $k + \delta k$ and truncates it to length $k$ in a way which uses the extra terms to correct the shorter expansion.  As the implementation described in \cite{CF} and provided in Chebfun \cite{ChebFun}
is an $O(\delta k^2)$ algorithm, it would seem that the CF method is not worth the computational effort either. In this work, however, we reduce this to $O(\delta k\ln \delta k)$, making its runtime similar to Chebyshev interpolation, and making its practical use worthwhile. These improvements in runtime can be seen in Figure \ref{fig:runtime2}.

In principle, it is possible to pre-compute minimax polynomials by choosing a fixed accuracy $\epsilon$, a fixed number of terms $k$, and then searching for the largest radius $r$ such that the $k$'th order minimax polynomial has maximum error less than $\epsilon$ on $[-r,r]$. For matrix polynomials, all eigenvalues of the matrix must then be on the interval $[-r,r]$ to ensure accuracy. This can be written as $\rho(\hat{H}) \leq r$ where $\rho(.)$ is the spectral radius, or the maximum absolute value of all eigenvalues of the matrix argument. This is also an approach used in \mono{scipy.sparse.linalg.expm\_multiply}, except using Maclaurin polynomials instead. Such an approach works well there, since the maximum stable $k$ is known, so any 
$\rho(\hat{H})$ that would require more terms must be split up in any case. However, for the CF method it is always possible to construct a problem with a greater $\rho(\hat{H})$ that needs more terms than have been calculated, so a method for generating polynomials on-demand is needed unless a similar splitting scheme is to be used.

As mentioned, the minimax problem requires a fixed interval, so an upper bound is needed on $\rho(\hat{H})$. In this case the induced $L_1$-norm is used as it can be calculated easily, and in practice has given a tighter bound than other viable options such as the Frobenius norm. Then, a polynomial $p$ is chosen 
such that $|p(x) - \exp(x)| < \epsilon$ for all $|x| < ||\hat{H}||_1$, $x \in \mathbb{R}$. The spectral radius $\rho(\hat{H})$ can be reduced slightly by shifting the energy spectrum of $\hat{H}_P$ by a constant at the start of the simulation so that the largest negative eigenvalue is equal in magnitude to the largest positive eigenvalue. In quantum mechanics, this changes the final state by a global phase factor, which doesn't affect any observables.

\begin{figure}
\includegraphics[width=0.90\linewidth]{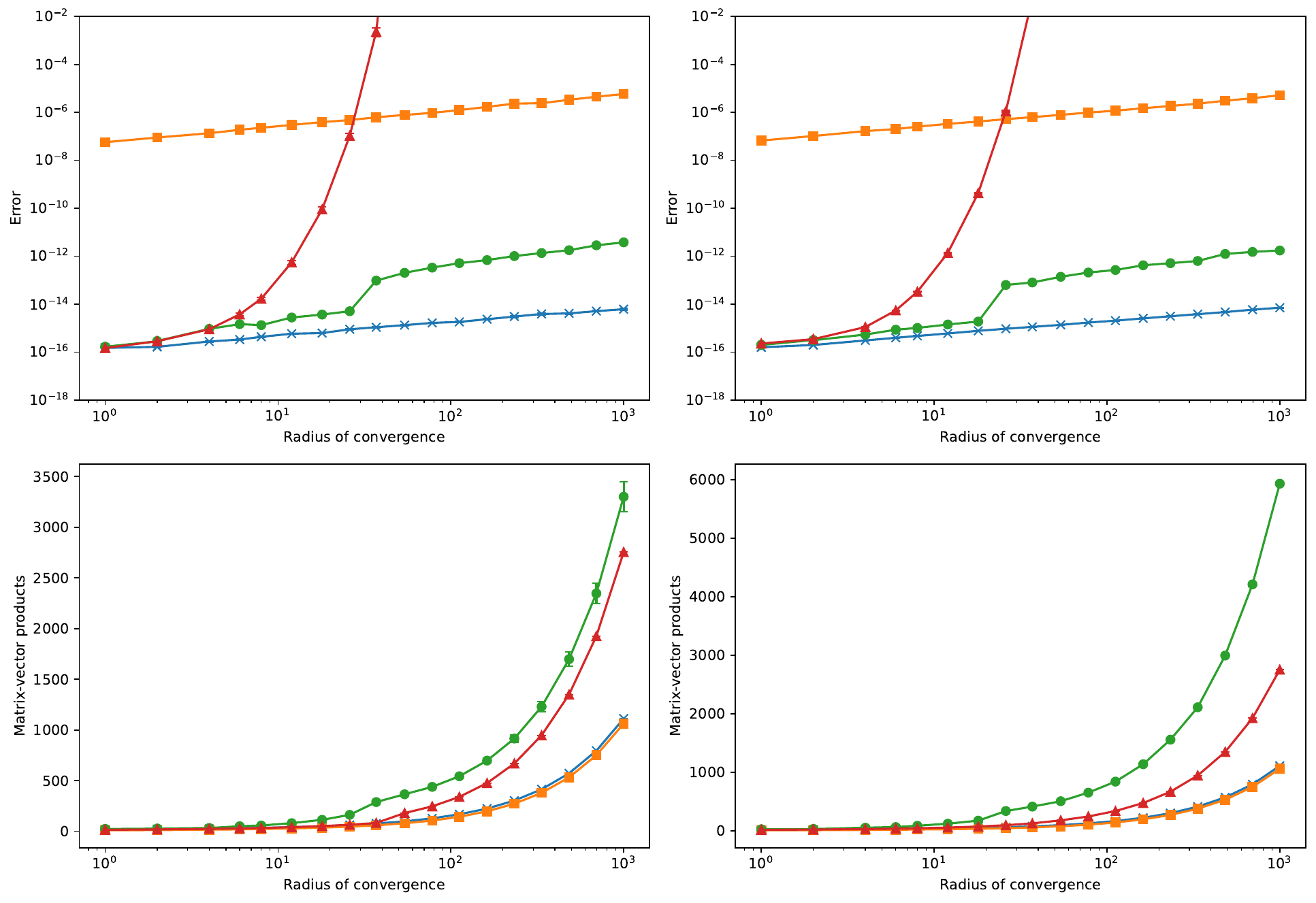}
\fcaption{A comparison of numerical errors and performance across different methods and problem instances for calculating $(\sin(H) + \cos(H)) \ket{\psi}$. The first column is tested on 100 5-qubit SK spin-glass problems from \cite{data}, where the parameters are sampled from a standard normal distribution. The second column is tested on 100 symmetric $32 \times 32$ matrices where all entries are sampled from a standard normal distribution, and then it is made symmetric by adding it to its own transpose. The first row shows how much numerical error is incurred, measured by calculating the euclidean distance of the output from an 80-bit precision baseline, as convergence radius increases. The second row shows how many matrix-vector products are needed, which is indicative of computational cost. Orange squares are the Clenshaw algorithm in single precision, red triangles are Taylor series in double precision, green circles are scipy's expm\_multiply function in double precision, and blue crosses are the Clenshaw algorithm in double precision. The Clenshaw algorithm shows excellent numerical stability while also providing a significant gain in performance.}
\label{fig:methods}
\end{figure}

\begin{figure}
\includegraphics[width=0.90\linewidth]{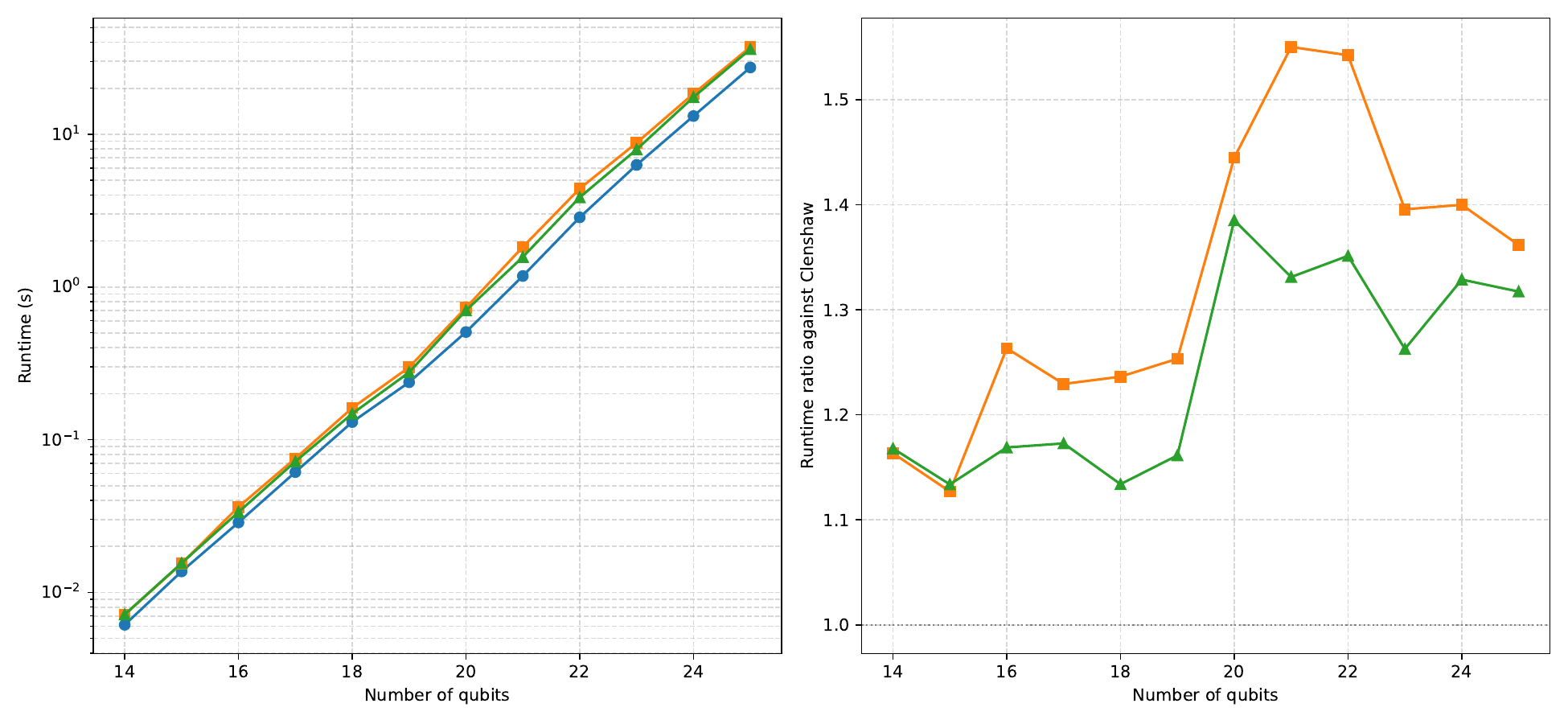}
\fcaption{The runtime needed to evaluate $p(\hat{H})$ for a set of 10 random SK spin-glass Ising Hamiltonians from \cite{data}, where parameters have been drawn from a standard normal distribution, and a degree-50 polynomial $p$ where the coefficients have been drawn from a standard normal distribution in the monomial basis. Orange squares show a naive Clenshaw implementation, green triangles show a naive Horner implementation, and blue circles show an optimised Clenshaw algorithm. On the right hand side, ratios of runtimes are shown for clarity. The optimised Clenshaw algorithm shows the most performance gain for large sizes, where it is most needed.}
\label{fig:runtime}
\end{figure}

For both the CF method and Chebyshev interpolation, the coefficients are returned in terms of the Chebyshev basis, which can be evaluated in an efficient and stable manner using the Clenshaw algorithm modified for matrix evaluations shown in Algorithm \ref{alg:clenshaw}.
\begin{algorithm}
\begin{minipage}{0.95\linewidth}
\begin{python}
def clenshaw(p, H, psi):
    """Evaluate p(H) @ psi using the Clenshaw algorithm, where p is
    expressed in the Chebyshev basis.
    Parameters:
    ----------
    p : list[float] of Chebyshev coefficients, where p[m] is the
    coefficient of T_m(x).
    psi : NDArray[np.floating] vector to which p(H) is applied.
    H : Any linear operator implementing `@`
    (matrix-vector product) and compatible with psi.
    """
    n = len(p)
    b_next = psi*0
    b_curr = p[n-1]*psi
    for r in range(n-2, 0, -1):
        b_curr, b_next = 2*H@b_curr - b_next + p[r]*psi, b_curr
    return H@b_curr - b_next + p[0]*psi
\end{python}
\fcaption{The Clenshaw algorithm for evaluating \mono{p(H)  @ psi}, in Python, with $p$ written in the Chebyshev basis.}
\label{alg:clenshaw}
\end{minipage}
\end{algorithm}
Figure \ref{fig:methods} shows that Chebsyhev-based evaluation converges faster and results in less numerical error than Taylor-based evaluation, and Figure \ref{fig:runtime} shows that the Clenshaw algorithm is not much slower than Horner's algorithm, although an optimised Horner's method would likely be slightly faster.

\subsection{Analysis of errors}
There are three main sources of error in polynomial approximation. These are: truncation error, coefficient error, and evaluation error. Let $f(x)$ be a smooth function on $[-1,1]$, which has an expansion in the form
\begin{align}
    f(x) = \sum_{m = 0}^\infty a_m \phi_m(x),
\end{align}
where $\phi_m$ are polynomials which form a basis. Here we are mainly concerned with the monomial basis $\phi_m(x) = x^m$ and the Chebyshev basis $\phi_m(x) = T_m(x)$. Let partial expansions be denoted by
\begin{align}
    p_k(x) = \sum_{m = 0}^k a_m \phi_m(x).
\end{align}
Truncation error is caused by only using a finite number of terms of an infinite expansion to approximate a function, and can be written as
\begin{align} \label{TruncError}
    |f(x) - p_k(x)| = \left|\sum_{m = k+1}^\infty a_m \phi_m(x)\right| \leq \sum_{m = k+1}^\infty |a_m\phi_m(x)|.
\end{align}
Truncation error is the easiest to control, as it can be reduced by using more terms, and can be estimated by calculating extra terms of the expansion and evaluating the upper bound in Eq.~\eqref{TruncError} directly.
Coefficient error is due to errors arising in the calculation of the polynomial coefficients. A routine for calculating $p_k$ will actually return a perturbed polynomial $p_k + \delta p_k$, which has coefficients $a_k + \delta_k$. The magnitudes of $\delta_k$ depends on the method used and the machine precision used. 
Evaluation error is due to errors that arise during the calculation of $p_k(x)$ at a specific value of $x$. This will be denoted $\tilde{p}_k(x)$. The total error can therefore be written as
\begin{align} \label{TotalError}
    (f(x) - \tilde{p}_k(x) - \tilde{\delta p}_k(x)) &= (f(x) - p_k(x)) + (p_k(x) - (p_k(x) + \delta p_k(x)) ) + (p_k(x) - \tilde{p}_k(x)) + (\delta p_k(x) - \tilde{\delta p}_k(x))\\
    &= \underbrace{(f(x) - p_k(x))}_{\color{red}{\text{Truncation}}} + \underbrace{\delta p_k(x)}_{\color{red}{\text{Coefficient}}} + \underbrace{(p_k(x) - \tilde{p}_k(x))}_{\color{red}{\text{Evaluation}}}+O(u^2),
\end{align}
where $u$ is the unit roundoff for the floating point type being used, which corresponds to a relative precision of $2^{-24}$ for single precision floating points.

\subsubsection{Analysis of coefficient error}
For a given function $f(x)$, an analytic form of the Taylor or Chebyshev expansion may be unknown or may be unreasonable to calculate. In this case, FFT based methods are the best choice for calculating the coefficients numerically. For the Taylor series, the substitution $x = e^{i \theta}$ transforms it into a Fourier series, making the application of an FFT clear. When an FFT based method is used to calculate $p_k$, coefficient errors can be large, and are difficult to bound. Table 5 in \cite{FFTerror} gives an upper bound on $\frac{|\delta p_k|}{M}$, where $M$ is the largest value of $|f(x)|$ over the Chebyshev nodes. For $N = 2^8$, the upper bound is $12310u$, with the worst known case being $1271u$, where $u = 2^{-53}$ for double precision floating point. Errors will typically be much better than the worst known case, but the bound represents the floor below which a coefficient is in principle indistinguishable from 0. Checks for convergence and error calculations using Eq \eqref{TruncError} must keep coefficient error in mind for this reason. As such, analytic expansions should be preferred when available.

Since the Chebyshev polynomials are an orthogonal basis, it is possible for some functions to evaluate the terms analytically by solving the resulting orthogonality integrals. Some functions for which analytic solutions are known are compiled in \cite{series}. Of interest here is Eq.~(5.18) in \cite{series}, which gives
\begin{align}
    e^{zx} &= I_0(z) + 2\sum_{m = 1}^\infty I_{m}(z) T_m(x) \\
    \implies e^{iax} &= J_0(a) +  2\sum_{m = 1}^\infty i^m J_{m}(a) T_m(x), \label{eq:cheb}
\end{align}
where $J_m$ is the Bessel function of the first kind, and $I_m$ is the modified Bessel function of the first kind. Both these series can be calculated efficiently with Miller's recurrence algorithm \cite{Miller} for the Bessel functions. The error analysis points out that during the quickly decaying part of the series, errors are suppressed by the same decay factor, meaning the resulting relative error on each coefficient is only what is introduced in that step, i.e., $\approx u$. This does not hold in the oscillatory section, however, where errors will build up as normal. This is still an improvement over the FFT, and having access to more accurate coefficients makes analysing truncation error much simpler. Since this step is comparatively cheap due to being purely scalar, it is recommended to use higher numerical precision to ensure all coefficients are correct to a relative precision of $\approx u$. This is achieved here with the x86 extended precision format. The error introduced into the final calculation due to coefficient error is then $|\delta p_k(x)| \lesssim ku $.

\subsubsection{Analysis of truncation error}
Truncation error can be controlled to be below any given tolerance $\epsilon$. Let $f(x) = \exp(irx)$ on the domain $[-1,1]$, for some $r > 0$. Let $P_k$ be the $k^{th}$ order Taylor series of $f$, and require that $P_k$ has truncation error below $\epsilon$. The error on the Taylor series is bounded by $\frac{r^{k+1}}{(k+1)!}$. The value of $k$ should be chosen such that the right hand side is smaller than the desired accuracy, $\epsilon$. To do this, set $\frac{r^{k+1}}{(k+1)!} = \epsilon$, then take the logarithm of both sides and apply Stirling's approximation to get an estimate.
\begin{align*}
    (k+1)\ln(r) - \ln((k+1)!) \approx (k+1)\ln(r) - (k+1)\ln(k+1) + (k+1) &= \ln(\epsilon)\\
    \implies (k+1)(\ln(r) - \ln(k+1) + 1) &= \ln(\epsilon) \\
    \implies \frac{k+1}{re}(\ln(\frac{k+1}{re})) &= -\frac{1}{re}(\ln(\epsilon)) \\
    \implies k &= re^{1 + W(\frac{-\ln(\epsilon)}{re})} - 1 \\
    \implies k & \sim re \text{ as } r \rightarrow \infty,
\end{align*}
where $W$ is the principal branch of the Lambert $W$ function.
This is also roughly the number of terms needed for a Chebyshev series to converge to an accuracy $\epsilon$ on the larger interval $[-2,2]$, as \cite{interpolation} gives the error term for Chebyshev interpolation as $\frac{r^{k+1}}{2^{k} k!} = \frac{2}{k!} \left(\frac{r}{2}\right)^{k+1}$. 

This is not a tight bound, however, as an asymptotic expansion by Debye gives
\begin{align}
    J_v(v \text{sech}(\alpha)) \sim \frac{e^{v(\tanh(\alpha) -\alpha)}}{\sqrt{2 \pi v \tanh(\alpha)}},
\end{align}
which for any $\alpha > 0$ tends to 0 exponentially as $v \rightarrow \infty$. Comparing with  Eq.~\eqref{eq:cheb}, this means that the $r(1 + \alpha)$'th term of $e^{irx}$ approaches 0 exponentially fast for any $\alpha > 0$ as $r \rightarrow \infty$. This means that, asymptotically, only $r + o(r)$ terms are needed, a factor of $e$ less than the Taylor series. This asymptote is approached rather slowly, as can be seen in Fig.~\ref{fig:convergence}. As a comparison, \cite{Higham} shows in Table 3.1 that 550 terms would be needed for convergence on the interval [-99,99] due to the need for scaling and splitting, whereas Chebyshev interpolation only requires 151 terms. Note that \cite{Higham} converts the truncation error to a backwards error and then requires that backwards error to be below a given tolerance. Their argument is entirely general and can be applied here too, but here a bound on the forwards error is used as it is only necessary that the truncation error is small compared to the other sources of error.

\subsubsection{Analysis of evaluation error} \label{sec:eval}
Evaluation error is the trickiest to analyse, and usually the most significant. Many different approaches have been taken to analyse the Clenshaw algorithm, both in the scalar and matrix cases \cite{bakshi2024,Smoktunowicz2002,smoktunowicz2013,generalbackwards}. However, the tightest known bound for the scalar case is given by \cite{interval}, building on ideas introduced in \cite{Elliott_1968}.  Here we extend it to the matrix case. A naive attempt at analysing Algorithm \ref{alg:clenshaw} will show that each term's error is roughly the sum of the errors of the previous two terms, and therefore grows exponentially like the Fibonacci sequence. This would clearly not be a stable algorithm, but empirical evidence such as Fig.~\ref{fig:methods} shows that the algorithm is stable in practice. 

It is first useful to write out Algorithm \ref{alg:clenshaw} in mathematical notation for the scalar case.
Let 
$$p_n(x) = \sum_{k = 0}^n c_k T_k(x)$$ be a truncated Chebyshev series. The polynomial $p_n(x)$ can be evaluated with the recurrence
\begin{align*}
    b_{n+1} &= 0, b_n = c_n,\\
    b_k &= 2xb_{k+1} - b_{k+2} + c_k,\\
    p_n(x) &= xb_1 - b_2 + c_0.
\end{align*}
As each $b_k(x)$ is calculated, some numerical error $\epsilon_k$ is introduced, and then carried forward throughout the calculation. The technique for analysing this in \cite{Elliott_1968} is to consider $\epsilon_k$ as coefficient error, that is to say that adding an error of $\epsilon_k$ at each stage leads to the same result as if the polynomial were defined as 
$$\tilde{p}_n(x) = \sum_{k = 0}^n (c_k + \epsilon_k) T_k(x),$$
and then the result was calculated with exact arithmetic.
The error is therefore $(\tilde{p}_n - p_n)(x)$. If an upper bound $e_k \geq |\epsilon_k|$ can be provided at every step, then the error can be bounded by evaluating the polynomial
$$\sum_{k=0}^ne_k|T_k(x)| \geq \sum_{k=0}^n \epsilon_kT_k(x).$$
This is essentially what Algorithm 2 of \cite{interval} does, but it applies $|T_k(x)| \leq 1$ to simplify the calculation and uses interval arithmetic to return an interval which the solution lies in. For the matrix-vector case, the Clenshaw algorithm can be written as 
\begin{align}
    \vec{b}_{n+1} &= \vec{0}, \vec{b}_n = c_n \vec{v}, \nonumber\\
    \vec{b}_k &= 2H\vec{b}_{k+1} - \vec{b}_{k+2} + c_k \vec{v}, \label{eq:recursion}\\
    p_n(H)\vec{v} &= H\vec{b}_1 - \vec{b}_2 + c_0 \vec{v}. \nonumber
\end{align}
Now, as each $\vec{b}_k$ is calculated, an error $\vec{{\epsilon}}_k$ is introduced. As this error is a vector, it cannot be absorbed into the coefficients as before, since they are still scalar.  However, consider using the Clenshaw algorithm to evaluate $T_m(H)\vec{\epsilon}_m$:
\begin{align*}
    \vec{b}_{m+1} &= \vec{0}, \vec{b}_m = \vec{\epsilon}_m,\\
    \vec{b}_k &= 2H\vec{b}_{k+1} - \vec{b}_{k+2} + \vec{0},\\
    p_n(H)\vec{\epsilon}_m &= H\vec{b}_1 - \vec{b}_2 + \vec{0}.
\end{align*}
By linearity, it can be seen that adding $T_m(H)\vec{\epsilon}_m$ to the final result is equivalent to introducing an error of $\vec{\epsilon}_m$ to $\vec{b}_m$ and then propagating it through the calculation, resulting in a total evaluation error of
$$\sum_{k=0}^n T_k(H)\vec{\epsilon}_k.$$
As in the scalar case, upper bounds can then be used. Some results and notation from \cite[Section~3.5]{HighamBook} will be used. Define $u$ as the unit roundoff, either $2^{-24}$ or $2^{-53}$ using single or double precision respectively, then define $\gamma_n = \frac{nu}{1 - nu}$. The error in the 2-norm on forming the matrix-vector product $H\vec{v}$ can be bounded as $||\tilde{H}\tilde{\vec{v}}- H \vec{v}||_2  \leq \sqrt{d}\gamma_d ||H||_2 ||\vec{v}||_2$, where $d$ is the maximum number of non-zero elements in a column of $H$. The error in forming the vector sum $\vec{v} + \vec{w}$ can be bounded by $|(\tilde{\vec{v}} + \tilde{\vec{w}}) - (\vec{v} + \vec{w})|_2 \leq u|\vec{v} + \vec{w}|_2$. Applying this to Eqs.~\eqref{eq:recursion} gives,
$$||\vec{\epsilon}_k||_2 \leq u||b_{k}||_2 + u||c_k \vec{v} - b_{k+2}|| + \sqrt{d} \gamma_d ||H||_2 ||\vec{b}_{k+1}||_2.$$
Since $H$ has been normalised to have $||H||_2 \leq 1,$ then
$$\sum_{k=0}^n T_k(H)\vec{\epsilon}_k \leq \sum_{k=0}^n ||T_k(H)\vec{\epsilon}_k||_2 \leq \sum_{k=0}^n ||T_k(H)||_2||\vec{\epsilon}_k||_2 \leq \sum_{k=0}^n ||\vec{\epsilon}_k||_2 \leq \sum_{k=0}^n u||b_k||_2 + u||c_k \vec{v} - b_{k+2}||+ \sqrt{d} \gamma_d ||\vec{b}_{k+1}||_2.$$

Note that the error will be slightly different depending on the order in which \eqref{eq:recursion} is evaluated, however, all three terms are expected to be similar in magnitude. Compared to Horner's method, which is analysed in \cite{Higham}, the evaluation error for the Clenshaw algorithm is essentially equivalent, scaling linearly with number of terms and being controlled by the size of the intermediate vectors that appear. A major advantage of the Clenshaw algorithm is that, for functions like $\exp(ix)$, the intermediate vectors $\vec{b}_k$ remain small, whereas in Horner's method the terms will blow up and then decay. Figure \ref{fig:methods} shows how errors compare empirically for Horner's method, Algorithm \ref{alg:horner} vs Clenshaw, Algorithm \ref{alg:clenshaw}.

\subsection{Improving on Chebyshev Expansion}
Approximating a function via polynomial interpolation at a fixed set of interpolation nodes will not, in general, be optimal. For a given set of nodes, how far the resulting polynomial can be from the minimax polynomial is described by the Lebesque constant, $\Lambda (D)$, where $D$ are the nodes. For the Chebyshev nodes, the resulting polynomial's accuracy is worse than the true optimum \cite{lebesgue} by a factor of up to $\Lambda(D) + 1 = \frac{2}{\pi} \ln(k+1) + 2$. To find the true optimum, two methods are presented. The first is the Remez algorithm \cite[Chapter 8]{Powell}. The goal of the Remez algorithm is to find a polynomial $p^*$ of degree $k$ such that 
$$\max_x |f(x) - p^*(x)| = ||f - p^*||_\infty$$ is minimised over the domain of $f$. It is known that $p^*$ satisfies the equi-oscillation theorem, and it is called the minimax polynomial of degree $k$. Suppose that $||f - p^*||_\infty = E$, the equi-oscillation theorem states that $f(x) - p^*(x)$ oscillates between $E$ and $-E$ repeatedly, with $|f(x) - p^*(x)|$ becoming equal to $E$ at $k + 2$ points \cite[Theorem 7.2]{Powell}. The Remez algorithm requires a starting polynomial, with the Chebyshev interpolant being the typical choice. Then, the $n+2$ extrema of the error curve are calculated. An interpolation problem is then set up, where $p(x_i)^* \pm E= f(x_i)$ at each extremum $x_i$, where the values $x_i$ are in increasing order. This is a well defined linear problem as there are $k+2$ extremum, and $k+2$ variables, as $E$ is also unknown. Solving this linear problem gives a new polynomial, and the process is then repeated until convergence. The implementation of the Remez algorithm used in our work uses a companion matrix approach to find the roots of the derivative of $p^{*}$, and then uses barycentric interpolation to solve the interpolation problem. Both these techniques are $O(n^2)$. Fast polynomial interpolation algorithms and fast multipoint evaluation algorithms exist that are $O(n \ln^2\/n)$, which could be leveraged to reduce the iterations to $O(n \ln^2\/n)$ \cite[Chapter 10]{MCAlg}. 

The second method for approximating $p^*$ is the Carath\'eodory-Fej\'er method \cite[Chapter 20]{ATAP}. Let $p_k$ and $p_{k + \delta k}$ be Chebyshev series of degree $k$ and $k + \delta k$ respectively. To correct the series of degree $k$, form a Hankel matrix using the last $\delta k$ coefficients of $p_{k + \delta k}$, and then find the largest magnitude eigenvalue, $\lambda$, and its eigenvector, $\vec{u}$. The eigenvalue estimates the minimax error, and the error curve is described by a Blaschke product which uses $\vec{u}$. Specifically, let $u(z)$ be the polynomial which has the entries of $\vec{u}$ as its coefficients, %
$$u(z) = \sum_{i = 0}^{\delta k} u_i z^{i}.$$ %
Then $b(z) = \frac{z^{\delta k} u(z^{-1})}{u(\bar{z})}$ is the relevant Blaschke product. The Laurent expansion of $b(z)$ then provides the corrections to $p_k(x)$ after mapping back to the unit interval with $x = \frac{1}{2}(z + z^{-1})$. Polynomial division can be done in $O(n \ln(n))$ runtime since division is pointwise after a Fourier transform.  However, there is a simpler way of achieving $O(n \ln(n))$ runtime which is faster in practice. An iterative procedure to find the Laurent expansion of $b(z)$ is  provided in \cite{CF},
\begin{align} \label{eq:iter}
    b_{m} = \frac{-1}{u_1}(b_{m+1} u_2 + b_{m+2} u_3 + \cdots + b_{m+\delta k - 1} u_{\delta k}).
\end{align}
On paper, this is $O(k^2)$ as it is $O(k)$ per iteration, and there are $k$ iterations. In practice, the entries of the Laurent expansion will often decay rapidly and allow an early exit from the iteration. The application of \eqref{eq:iter} can be reinterpreted as a matrix multiplication, with repeated iterations being equivalent to the power method for finding the leading eigenvector. So long as the matrix has spectral radius less than 1, then the iterations decay exponentially, and early exit will be achieved in a logarithmic number of iterations.

The runtimes of a set of approximation methods are compared in Figure \ref{fig:runtime2}.

\subsection{Numerical implementation methods}
The simulations presented in Figs.~\ref{fig:convergence}--\ref{fig:runtime} were run on a machine with a Ryzen 9 7950X CPU, an RTX 4090 GPU (with 24GB of VRAM) and 64GB of 6000MHz DDR5 RAM. For simulations on CPU, there is 1MB of L2 cache per core which can handle the cases up to $n = 16$, and 64MB of global L3 cache which can handle cases up to $n = 22$, after which cache misses are inevitable. On this machine, the VRAM has a peak bandwidth of roughly 20 times the main memory.

Profiling of the code shows that the algorithm becomes bottlenecked by memory access as the number of qubits grows, so for optimization, the focus has been on RAM usage and cache-hit behaviour. For the SK spin glasses in \cite{data}, a naive implementation would use $O(nN)$ RAM since $\hat{H}$ has that many nonzero elements, but it is possible to reduce RAM usage to $O(N)$ by using a technique similar to the fast Walsh-Hadamard Transform to multiply $\ket{\psi}$ by $\hat{H}_G$ without explicitly writing down $\hat{H}_G$. This follows from a recursive definition of $\hat{H}_G$,
$$ \hat{H}_G(n) = \begin{cases}
I_2\otimes \hat{H}_G(n-1)+\hat{H}_G(1)\otimes I_{2^{n-1}} & \text{if } n>1\\
\begin{bmatrix}
0 & 1\\
1 & 0
\end{bmatrix}
&\text{if }n=1
\end{cases}
$$
which allows a divide-and-conquer strategy to be used. This allows $\hat{H}\ket{\psi}$ to be evaluated as $\hat{H}_I\ket{\psi} + f(\ket{\psi})$, where $f$ is the transform and $\hat{H}_I$ is diagonal.

In the naive implementation of the Clenshaw algorithm in \ref{alg:clenshaw}, it first appears that an extra $16N$ bytes of RAM are needed compared to \ref{alg:horner}, corresponding to the \mono{b1} and \mono{b2} variables. $8N$ bytes can be saved by not using extra working space for \mono{H @ b1}, and instead directly calculating \mono{b2 += H @ b1}. In principle, changing the polynomial basis to use Horner's method with the minimax polynomial will be faster, but this reintroduces the numerical instabilities of the Taylor series approach.

This means that 4 vectors need to be stored for the Clenshaw algorithm, three of which are complex values. These are the diagonal of $\hat{H}_I$, the state vector $\ket{\psi}$, and the vectors \mono{b1} and \mono{b2} used for the Clenshaw algorithm. There are $N$ entries of each, with each entry requiring $4$ bytes in single-precision floating point, or $8$ bytes for a complex vector. This adds up to needing $28N$ bytes of RAM, or $56N$ bytes in double precision. On both CPU and GPU, memory accesses are the performance bottleneck, so reducing the RAM needed allows more data to be held in cache and reduces accesses to main memory.

The current implementation of \mono{H @ b1} uses explicit vectorisation to make use of the permute instructions. On a CPU with AVX512 available, these instructions can effectively perform the Walsh-Hadamard transform on 16 entries in a single instruction. The current code is designed so that while these slow permute instructions are executing, out-of-cache data is loaded so the penalty from cache misses is mitigated.
The implementation here also keeps $\hat{H}$ as a real data type, so each $\hat{H}v$ evaluation uses half as many operations as passing $i\hat{H}$ into a complex exponential function does.

\section{Discussion and further work}\label{sec:diss}

We have presented and discussed improvements and optimizations for evaluating the action of a matrix exponential on a vector.  The associated C++ code \cite{AsaCode} provides significant improvements over commonly available open-source implementations, and algorithmic stability has now been proven. Figure \ref{fig:methods} shows a significant advantage for the matrix Clenshaw algorithm over Horner's algorithm both in performance and numerical stability, and so should be preferred when available. On the other hand, Figure \ref{fig:convergence} shows a rather unconvincing argument for using true minimax polynomials against Chebyshev series. Nevertheless, it is worth keeping in mind that using the true minimax polynomial still reduces truncation error significantly. For calculations that involve chaining together many such operations (e.g., multi-stage quantum walks, or the extension to time-dependent Hamiltonians below), or approximations of functions that aren't infinitely smooth, the reduced truncation error will be more noticeable. In situations where polynomials are precomputed, or their computation time is negligible, true minimax polynomials can provide benefits with few drawbacks. The technique used to analyse evaluation error in Section \ref{sec:eval} can be applied to other three-term recurrences, such as those used for other common polynomial bases as in \cite{Laguerre}.
There are several useful directions for extending the current work.  We discuss extra considerations for GPU implementation, and extensions to time dependent Hamiltonians below.  Enhancements to the associated C++ code \cite{AsaCode} could include adding Python bindings, and making it easier to swap between CPU and GPU based computation and 32 or 64-bit precision.

\subsection{GPU implementation}

\begin{figure}
\includegraphics[width=0.90\linewidth]{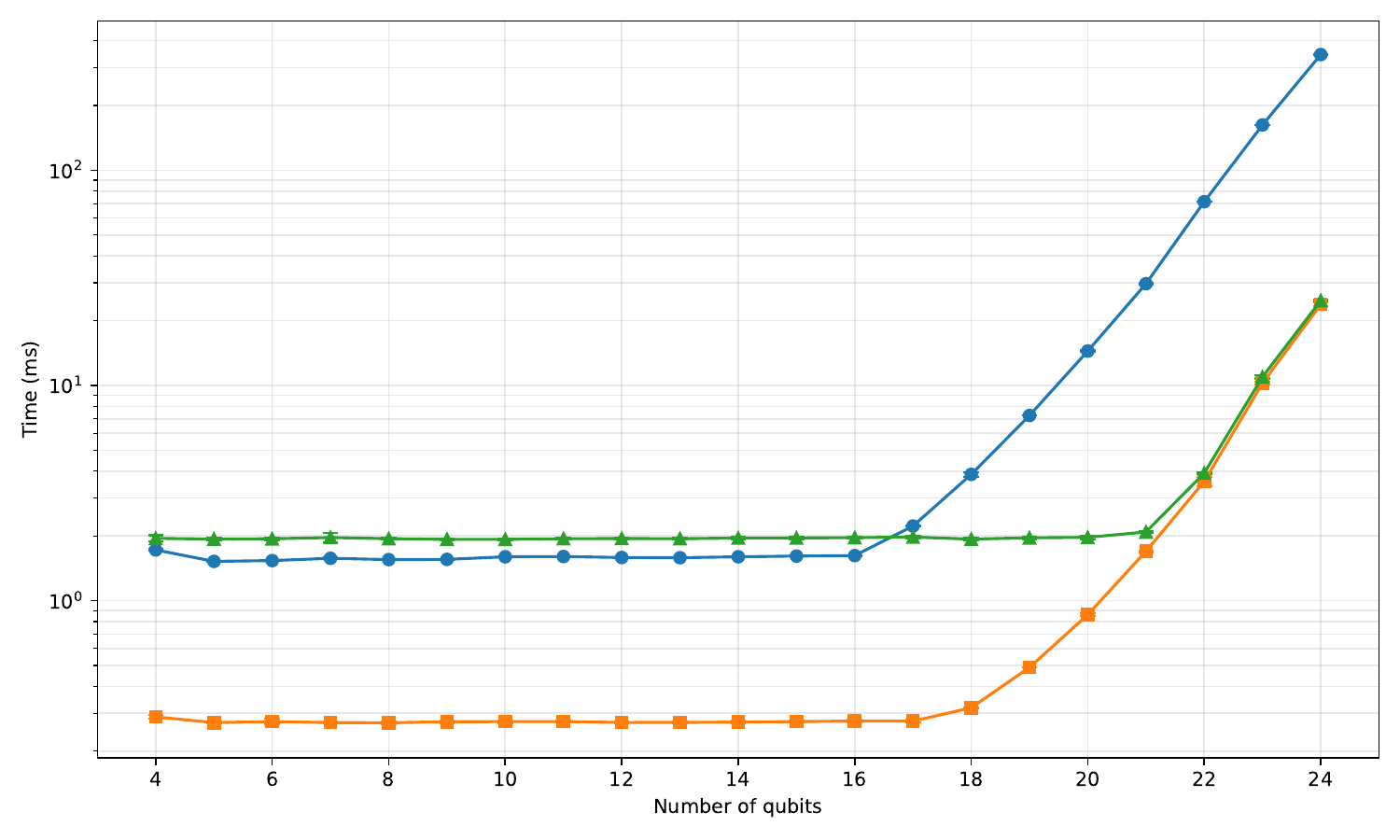}
\fcaption{The runtime needed to evaluate $p(\hat{H}) \ket{\psi}$ on a GPU for a set of 20 degree-50 polynomials $p$ where the coefficients have been drawn from a standard normal distribution in the monomial basis. Blue circles show sparse matrix multiplication with CuPy, orange squares show a matrix-free approach using a CuPy RawKernel, and green triangles show a matrix-free approach using a Taichi kernel. In the large $n$ limit, Taichi's performance matches CUDA, but for moderate $n$ the kernel launch overhead is too significant.}
\label{fig:GPU}
\end{figure}

GPU implementation has been tested, but poses some difficulties. The library CuPy \cite{CuPy} provides drop-in replacements for many NumPy and SciPy functions that work on GPU, including sparse linear algebra. This can give significant speedup for sparse operations, however, taking advantage of the matrix-free approach becomes difficult. CuPy does have the \mono{RawKernel} function, which allows writing of CUDA code directly. However, code written this way will only compile for Nvidia GPUs. An alternative with much wider hardware support is Taichi \cite{Taichi}, 
which provides a way of writing massively parallel code which can be compiled for a wide range of hardware with the available backends. Performance is similar to using CUDA directly, however, the biggest downside of Taichi is that kernels take much longer to launch, the effect of which can be seen in Figure \ref{fig:GPU}. When originally tested for this project, kernel launch times made GPU acceleration appear slower for the values of $n$ being examined, and so CPU implementation was favoured. For large $n$, kernel launch times are negligible and so the expected benefit from GPU acceleration can be achieved.

To overcome the issue of kernel launch times, it is necessary to batch together as much work as possible in each kernel call. One way of doing this for our $\hat{H}$ is to take $2^{24 - n}$ problem instances at once, such that the concatenated state vectors total $2^{24}$ entries, and then to operate on this concatenated state. This requires some care, as each problem instance has its own polynomial associated with it, and its own $\hat{H}_P$ matrix. Since $\hat{H}_P$ is diagonal, they can also be concatenated in the same way. As for the different polynomials, the suggested approach is to find the problem in the batch that requires the polynomial with the largest radius of convergence, and to use that polynomial for every problem instance in the batch. If one problem has much larger spectral radius than the others, this can cause an unnecessary amount of work to be done, so it may be worth passing through all problems once and intentionally batching together problems with similar spectral radius.

In order to reduce the amount of data sent between the main memory and GPU, it would also be possible to only send the coupling parameters of the Ising problems to the GPU, and then the calculation of the entries of $\hat{H}_P$ could be done entirely on the GPU.

\subsection{Extension to time dependent Hamiltonians}
\begin{figure}
\includegraphics[width=0.90\linewidth]{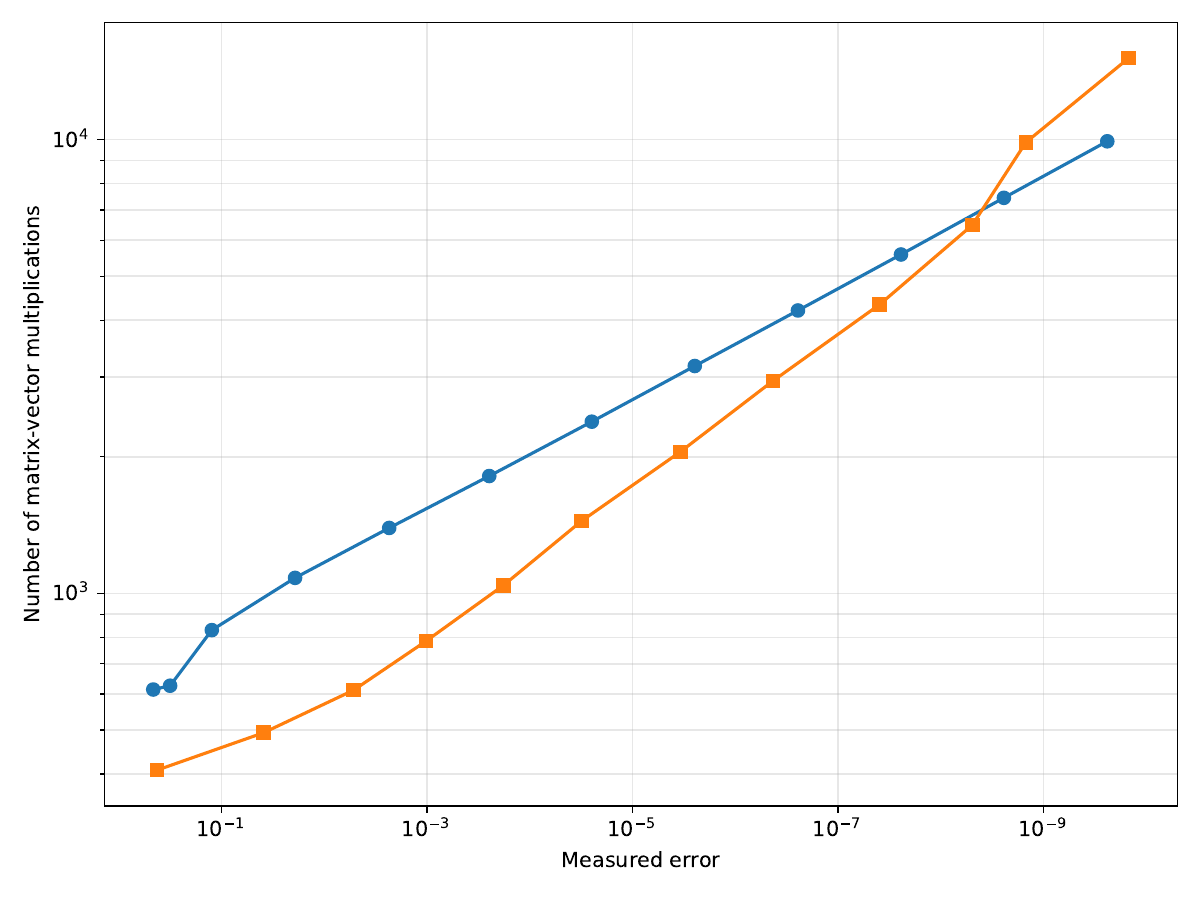}
\fcaption{The number of matrix-vector multiplications needed to simulate a quantum anneal for an 8-qubit SK spin-glass Hamiltonian from \cite{data} at different levels of measured error. Blue circles show the DOP853 method provided by scipy, and orange squares show the commutator-free method with Chebyshev polynomials. It is currently unclear why performance degrades suddenly for errors below $10^{-8}$}
\label{fig:CFET}
\end{figure}

There are various techniques for approximating the evolution of a quantum system under a time-dependent Hamiltonian using time-independent steps. A widely used example is Suzuki-Trotter decomposition, however, convergence is poor so many stages are needed. A method is described in \cite{commutator} which is built on the Magnus expansion, but manages to eliminate any commutators from explicitly appearing in the exponential steps. A weakness of this method is that choosing timesteps can be an expensive process. The authors suggest a method that requires two runs at large timesteps along with known error scaling to estimate the timestep needed for a given error tolerance, but this can increase running time considerably. Some significant improvements to this algorithm are possible in the case of quantum annealing in particular, due to bounds on the commutators that appear in the error formulae being possible to calculate in $O(\text{poly}(n))$ time. For the fourth order method in \cite{commutator}, an explicit error term is provided, and replacing the commutators that appear with upper bounds allows the error term to be calculated in $O(1)$ time. The constants provided are the result of an optimisation process, but it is viable to optimise the constants for the specific problem at hand. Secondly, being able to evaluate the error term directly in $O(1)$ time allows for choosing timesteps without running the algorithm at all. These changes make the method compare very favourably against more traditional ODE methods, and preliminary results in Figure \ref{fig:CFET} show this technique requiring much fewer matrix-vector multiplications of the ODE-based methods provided in QuTip \cite{QuTiP} at moderate accuracy. Measuring matrix-vector products also ignores the significant amount of vector-vector operations that methods like DOP853 perform, so real-time performance compares even more favourably.

\begin{acknowledgments}
AH is funded by UKRI EPSRC PhD studentship number 2745408. 
VK funded by UKRI EPSRC grants EP/T026715/2, 
EP/Z53318X/1, and the UKRI DRI and STFC funded CCP-QC Bridge Project -- extended case studies for neutral atom hardware.
\end{acknowledgments}


\bibliography{apssamp}

\end{document}